\documentclass[11pt,leqno,a4paper]{article}

\long\def\forget#1{}

\forget{
\pdfoutput=1

\usepackage[pdftex,bookmarks,bookmarksopen,bookmarksnumbered]{hyperref}
\hypersetup{%
 pdftitle     = {Change of Coefficients for Drinfeld Modules and Abelian Sheaves}
 pdfauthor    = {Urs Hartl, Markus Hendler},
 pdfcreator = {pdfLaTeX},
 pdfproducer = {pdfLaTeX with hyperref}
}
}

\forget{
\def\theenumi{(\alph{enumi})}

\def\p@enumii{\theenumi}
}

\usepackage{amsmath}
\usepackage{amssymb}
\usepackage{amscd}
\usepackage{textcomp}
\usepackage{ifthen}

\newcounter{zahl}

\newcommand{\DS}{\displaystyle}

\DeclareMathOperator{\Aut}{Aut}

\DeclareMathOperator{\End}{End}

\DeclareMathOperator{\Frob}{Frob}

\DeclareMathOperator{\Hom}{Hom}
\newcommand{\CHom}{{\cal H}om}

\DeclareMathOperator{\Isom}{Isom}

\DeclareMathOperator{\Spec}{Spec}

\DeclareMathOperator{\Sym}{Sym}

\newcommand{\alg}{{\rm alg}}

\DeclareMathOperator{\coker}{coker}

\DeclareMathOperator{\id}{\,id}

\DeclareMathOperator{\ord}{ord}

\DeclareMathOperator{\Sch}{\CS \!{\it ch}}
\DeclareMathOperator{\Sets}{\CS \!{\it ets}}
\newcommand{\DrMod}{\CalD r\text{-}A\text{-}\CM od\,}
\newcommand{\DrModp}{\CalD r\text{-}A'\text{-}\CM od\,}
\newcommand{\AbSh}{C\text{-}\CA b\text{-}\CS h}
\newcommand{\AbShp}{C'\text{-}\CA b\text{-}\CS h}
\newcommand{\DrSht}{C\text{-}\CalD\CA\text{-}\CS ht\,}
\newcommand{\DrShtp}{C'\text{-}\CalD\CA\text{-}\CS ht\,}

\newcommand{\DRM}{\mbox{\rm\ul{Dr-$A$-Mod}}}
\newcommand{\DRMp}{\mbox{\rm\ul{Dr-$A'$-Mod}}}
\newcommand{\CAB}{\mbox{\rm\ul{$C$-Ab-Sh}}}
\newcommand{\CpAB}{\mbox{\rm\ul{$C'$-Ab-Sh}}}
\newcommand{\DAS}{\mbox{\rm\ul{$C$-DA-Sht}}}
\newcommand{\DASp}{\mbox{\rm\ul{$C'$-DA-Sht}}}

\renewcommand{\phi}{\varphi}
\renewcommand{\epsilon}{\varepsilon}

\usepackage{amsfonts}
\newcommand{\BOne}{\underline{\hbox{\rm1\kern-2.5pt l\kern.9pt}}}
\newcommand{\OOne}{\BO\kern-7.7pt\raisebox{1.3pt}{\underline{\phantom{{$\BO$}}}}}

\newcommand{\BF}{{\mathbb{F}}}
\newcommand{\BG}{{\mathbb{G}}}

\newcommand{\BO}{{\mathbb{O}}}
\newcommand{\BP}{{\mathbb{P}}}

\newcommand{\BZ}{{\mathbb{Z}}}

\newcommand{\CA}{{\cal{A}}}

\newcommand{\CalD}{{\cal{D}}}
\newcommand{\CE}{{\cal{E}}}
\newcommand{\CF}{{\cal{F}}}
\newcommand{\CG}{{\cal{G}}}
\newcommand{\CH}{{\cal{H}}}

\newcommand{\CJ}{{\cal{J}}}

\newcommand{\CL}{{\cal{L}}}
\newcommand{\CM}{{\cal{M}}}

\newcommand{\CO}{{\cal{O}}}

\newcommand{\CS}{{\cal{S}}}
\newcommand{\CT}{{\cal{T}}}

\def\longto{\longrightarrow}

\def\isoto{\arrover{\sim}}

\newbox\mybox
\def\arrover#1{\mathrel{
       \setbox\mybox=\hbox spread 1.4em{\hfil$\scriptstyle#1$\hfil}
       \vbox{\offinterlineskip\copy\mybox
             \hbox to\wd\mybox{\rightarrowfill}}}}

\usepackage[curve,matrix,arrow,cmtip]{xy}

\newdir^{ (}{{}*!/-3pt/\dir^{(}}    
\newdir^{  }{{}*!/-3pt/\dir^{}}    
\newdir_{ (}{{}*!/-3pt/\dir_{(}}    
\newdir_{  }{{}*!/-3pt/\dir_{}}    

\let\setminus\smallsetminus

\newcommand{\es}{\enspace}

\newcommand{\dual}{^{^\vee}}

\newcommand{\mal}{^{^\times}}

\newcommand{\invlim}[1][]{\ifthenelse{\equal{#1}{}}
{\DS \lim_{\longleftarrow}}
{\DS \lim_{\underset{#1}{\longleftarrow}}}
}
\newcommand{\dirlim}[1][]{\ifthenelse{\equal{#1}{}}
{\DS \lim_{\longrightarrow}}
{\DS \lim_{\underset{#1}{\longrightarrow}}}
}
\newcommand{\ul}[1]{{\underline{#1}}}

\newcommand{\wt}[1]{{\widetilde{#1}}}

\usepackage{amsthm}
\theoremstyle{plain}
\newtheorem{theorem}{Theorem}[section]

\newtheorem{lemma}[theorem]{Lemma}

\newtheorem{proposition}[theorem]{Proposition}

\theoremstyle{definition}
\newtheorem{definition}[theorem]{Definition}

\def\?{\ ???\ \immediate\write16{}%
\immediate\write16{Warning: There was still a question mark . . . }%
\immediate\write16{}}

\newcommand{\UCF}{\underline{\CF}\mathstrut}

\def\longto{\longrightarrow}

\begin{document}


\author{Urs Hartl, Markus Hendler
\footnote{Both authors acknowledge support of the Deutsche Forschungsgemeinschaft in form of DFG-grant HA3006/2-1}
}

\title{Change of Coefficients for Drinfeld Modules, \\ Shtuka, and Abelian Sheaves}


\maketitle

\begin{abstract}
\noindent
We study the passage from Drinfeld-$A'$-modules to Drinfeld-$A$-modules for a given finite flat inclusion $A\subset A'$. We show that this defines a morphism from the moduli space of Drinfeld-$A'$-modules to the moduli space of Drinfeld-$A$-modules which is proper but in general not representable. For Drinfeld-Anderson shtuka and abelian sheaves instead of Drinfeld modules we obtain the same results.

\end{abstract}

 
\section*{Introduction}

Throughout this article let $\BF_q$ be a finite field with $q$ elements and characteristic $p$ and let $C$ and $C'$ be two smooth projective geometrically irreducible curves over $\BF_q$. Let $\pi:C'\to C$ be a fixed finite morphism of degree $n$. Let $\infty\in C$ be a closed point which does not split in $C'$, that is, there is exactly one point $\infty'\in C'$ above $\infty$. Set $A:=\Gamma(C\setminus\{\infty\},\CO_C)$ and $A':=\Gamma(C'\setminus\{\infty'\},\CO_{C'})$, then $A'$ is a flat $A$-algebra via $\pi^\ast:A\to A'$. 

In this situation $\pi$ defines a \emph{restriction of coefficients functor} from Drinfeld-$A'$-modules over $S$ to Drinfeld-$A$-modules over $S$. This functor induces a morphism between the moduli spaces (moduli functors, or more sophisticated, moduli stacks) classifying Drinfeld-$A'$-modules, respectively Drinfeld-$A$-modules. We show in this article that this morphism is proper but not necessarily representable. Likewise we study the effect of $\pi$ on Drinfeld-Anderson shtuka, see Definition~\ref{Def5.1}, and on abelian sheaves, a notion introduced by the first author~\cite{HartlAbSh} as a higher dimensional generalization of Drinfeld modules, see Definition~\ref{Def1.5}. For the case of Drinfeld-Anderson shtuka we may even relax the condition on $\pi$ and drop the assumption on the ramification of $\infty$. The pushforward of sheaves along $\pi\times\id_S:C'_S\to C_S$ defines a \emph{restriction of coefficients functor} from Drinfeld-Anderson shtuka on $C'$ over $S$ to Drinfeld-Anderson shtuka on $C$ over $S$, respectively from abelian sheaves on $C'$ over $S$ to abelian sheaves on $C$ over $S$. Again this yields proper but in general not representable morphisms between the moduli spaces classifying Drinfeld-Anderson shtuka on $C'$, respectively Drinfeld-Anderson shtuka on $C$ and similarly for abelian sheaves.

Of course the results for Drinfeld modules, Drinfeld-Anderson shtuka, and abelian sheaves are strongly related by the fact that the category of Drinfeld-$A$-modules over $S$ is anti-equivalent to a full subcategory of the category of Drinfeld-Anderson shtuka on $C$ over $S$ and anti-equivalent to a full subcategory of the category of abelian sheaves on $C$ over $S$. Nevertheless we give proofs also for the case of Drinfeld modules since these are particularly simple. After recalling the definitions and some basic properties in Section~\ref{Sect1} we prove in Sections~\ref{Sect2}, \ref{Sect3}, and \ref{Sect4} the properness and non-representability results for Drinfeld modules, abelian sheaves, respectively Drinfeld-Anderson shtuka.

This article has its origin in a conversation with F.\ Breuer who mentioned to us a special case of the proof for properness in the case of Drinfeld modules. Our proof of Proposition~\ref{Prop2.3} below is a generalization of his. We like to express our gratitude to him.


\section{Drinfeld Modules, Shtuka, and Abelian Sheaves} \label{Sect1}

We retain the notation from the introduction. In addition, we set $\deg(\infty):=[\kappa(\infty):\BF_q]$ and we denote by $\ord_\infty$ the normalized valuation on the fraction field of $A$ associated with the place $\infty$. For an $\BF_q$-scheme $S$ we set $C_S:=C\times_{\BF_q}S$. Unless mentioned explicitly we make no noetherian assumption on $S$.

For an $\BF_q$-algebra $B$ we denote by $B\{\tau\}$ the non-commutative polynomial ring in the variable $\tau$ over $B$ with the commutation rule $\tau\,b=b^q\,\tau$ for all $b\in B$. As in \cite[\S 1]{Matzat} one sees

\begin{proposition}\label{Prop1.2}
There is an isomorphism of rings between $B\{\tau\}$ and $\End_{B,\BF_q}(\BG_{a,B})$ the ring of $\BF_q$-linear endomorphisms of the additive group scheme over $\Spec B$ given by mapping $\tau$ to the $q$-th power Frobenius of $\BG_{a,B}$. \qed
\end{proposition}

\begin{definition}\label{Def1.1} (Drinfeld~\cite[\S 5.B]{Drinfeld})\\ 
Let $S$ be an $\BF_q$-scheme and assume there is a morphism $c:S\to \Spec A$. Let $r$ be a positive integer. A \emph{Drinfeld-$A$-module of rank $r$ and characteristic $c$} over $S$ is a pair $(E,\phi)$ where $E$ is a commutative group scheme over $S$ and 
\[
\phi:\es A\longto\End_S(E)
\]
is a ring homomorphism from $A$ to the ring $\End_S(E)$ of endomorphisms of the $S$-group scheme $E$ such that
\begin{enumerate}
\item 
$E$ is Zariski locally on $S$ isomorphic to the additive group scheme $\BG_{a,S}$,
\item 
if $U=\Spec B$ is an affine open subset of $S$ and $\psi:E_U\isoto\BG_{a,U}$ is an isomorphism of $S$-group schemes then for each $a\in A\setminus\{0\}$
\[
\psi\circ\phi(a)\circ\psi^{-1}\es=\es\sum_{i=0}^{<\infty}\delta_i(a)\tau^i\es\in\es B\{\tau\}
\]
with $\delta_0(a)=c^\ast(a)$, $\delta_i(a)\in B\mal$ for $i=d(a):=-r \ord_\infty(a)\deg(\infty)$, and $\delta_i(a)$ nilpotent for $i>d(a)$.
\end{enumerate}
A \emph{morphism} of Drinfeld-$A$-modules $\epsilon:(E,\phi)\to (\wt E,\wt\phi)$ is a morphism of $S$-group schemes $\epsilon:E\to\wt E$ which satisfies $\wt\phi(a)\circ\epsilon=\epsilon\circ\phi(a)$ for all $a\in A$.
\end{definition}

If $f:S'\to S$ is a morphism of $\BF_q$-schemes we can pull back Drinfeld-$A$-modules $(E,\phi)$ over $S$ to Drinfeld-$A$-modules $(f^\ast E,f^\ast\phi)$ over $S'$.

The following proposition is due to Drinfeld~\cite[Propositions 5.1 and 5.2]{Drinfeld}

\begin{proposition}\label{Prop1.3}
Let $(E,\phi)$ be a Drinfeld-$A$-module of rank $r$ over $S$. Then Zariski locally on $S$ there exists an isomorphism $\epsilon:(E,\phi)\isoto(\BG_{a,S},\psi)$ of Drinfeld-$A$-modules where $\psi$ is of the \emph{standard form}
\[
\psi:\es A\longto\CO_S\{\tau\}\,,\quad\psi(a)\es=\es\sum_{i=0}^{d(a)}\delta_i(a)\tau^i
\]
with $d(a):=-r \ord_\infty(a)\deg(\infty)$ and $\delta_{d(a)}\in\CO_S\mal$. Moreover if $\psi(a)$ is of the described form for one $a\in A\setminus\BF_q$ then it already is for any $a\in A$. \qed
\end{proposition}

\begin{proposition}\label{Prop1.4}
The morphism $\pi:C'\to C$ defines a \emph{restriction of coefficients functor} $\pi_\ast:(E',\phi')\mapsto (E',\phi'\circ\pi^\ast)$ from Drinfeld-$A'$-modules of rank $r'$ over $S$ to Drinfeld-$A$-modules of rank $nr'$ over $S$, where $n$ is the degree of $\pi$.
\end{proposition}

\begin{proof}
The change of rank results from the fact that $n\ord_\infty(a)\deg(\infty)=\ord_{\infty'}(a)\deg(\infty')$ for all $a\in A$ since $\pi^{-1}(\infty)=\{\infty'\}$. The rest is clear from the definition.
\end{proof}

\noindent
{\it Remark.} 
Consider the moduli problem, that is, the contravariant functor
\begin{eqnarray*}
\DRM^r:\es\Sch_{/\Spec A}&\longto&\Sets\\[2mm]
(c:S\to\Spec A) &\mapsto &\Bigl\{\,\text{Isomorphism classes of Drinfeld-$A$-modules} \\[-2mm]
& & \text{\qquad of rank $r$ and characteristic $c$ over $S$}\,\Bigr\}
\end{eqnarray*}
from the category of schemes over $\Spec A$ to the category of sets. This functor is not representable (without adding level structures). Nevertheless the restriction of coefficients functor defines a \emph{restriction of coefficients morphism} 
\[
\pi_\ast:\es\DRMp^{r'}\es\longto\es\DRM^{nr'}\,,\quad(E',\phi')\es\mapsto\es\pi_\ast(E',\phi')\,.
\]

\medskip

\noindent
{\it Remark.} 
If we let $S$ vary, the category of Drinfeld-$A$-modules of rank $r$ becomes a stack $\DrMod^r$ for the fppf topology on the category of $\BF_q$-schemes. It is an algebraic stack in the sense of Deligne-Mumford~\cite{DM}, see Laumon~\cite[Corollary 1.4.3]{Laumon}. The restriction of coefficients functor defines a \emph{restriction of coefficients 1-morphism} $\pi_\ast:\DrModp^{r'}\to\DrMod^{nr'}$.

\bigskip

Next we study the analogous situation for abelian sheaves. This notion was introduced in \cite{HartlAbSh}. While Drinfeld modules are analogues for elliptic curves in the arithmetic of function fields, abelian sheaves are the appropriate analogues for abelian varieties as the results of \cite{HartlAbSh,BoHa} amply demonstrate.

Let $r$ and $d$ be positive integers and write $\frac{d}{r\deg(\infty)}=\frac{k}{\ell}$ with relatively prime positive integers $k$ and $\ell$.
Let $S$ be an $\BF_q$-scheme and fix a morphism $c:S\to C$. Let $\CJ$ be the ideal sheaf on $C_S$ of the graph of $c$. We let $\sigma:=\id_C\times\Frob_q$ be the endomorphism of $C_S$ that acts as the identity on the underlying topological space and on the coordinates of $C$ and as $b\mapsto b^q$ on the elements $b\in \CO_S$. Let $pr:C_S\to S$ be the projection onto the second factor. For an integer $m$ denote by $\CO_{C_S}(m\cdot\infty)$ the invertible sheaf on $C_S$ associated with the divisor $m\cdot\infty$ and set $\CF(m\cdot\infty):=\CF\otimes_{\CO_{C_S}}\CO_{C_S}(m\cdot\infty)$ for any sheaf of $\CO_{C_S}$-modules on $C_S$.

\begin{definition}\label{Def1.5}
An \emph{abelian sheaf $\ul{\CF}=(\CF_i,\Pi_i,\tau_i)$ on $C$ of rank $r$, dimension $d$, and characteristic $c$ over $S$\/} is a ladder of locally free sheaves $\CF_i$ on $C_S$ of rank $r$ and injective homomorphisms $\Pi_i$, $\tau_i$ of $\CO_{C_S}$-modules ($i\in \BZ$) of the form
\[
\begin{CD}
\cdots & @>>> & \CF_{i-1} & @>{\Pi_{i-1}}>> & \CF_i & @>{\Pi_i}>> & \CF_{i+1} & @>{\Pi_{i+1}}>> & \cdots \\
& & & & @AA{\tau_{i-2}}A & & @AA{\tau_{i-1}}A & & @AA{\tau_i}A & & \\
\cdots & @>>> & \sigma^\ast\CF_{i-2} & @>{\sigma^\ast\Pi_{i-2}}>> & \sigma^\ast\CF_{i-1} & @>{\sigma^\ast\Pi_{i-1}}>> & \sigma^\ast\CF_i & @>{\sigma^\ast\Pi_i}>> & \cdots 
\end{CD}
\]
subject to the following conditions (for all $i\in \BZ$):
\begin{enumerate}
\item \label{DefAbelianSheafCond1}
the above diagram is commutative,
\item \label{DefAbelianSheafCond2}
the morphism $\Pi_{i+\ell-1}\circ\ldots\circ\Pi_i$ identifies $\CF_i$ with the subsheaf $\CF_{i+\ell}(-k\cdot\infty)$ of $\CF_{i+\ell}$,
\item\label{DefAbelianSheafCond3}
$pr_\ast\coker\Pi_i$ is a locally free $\CO_S$-module of rank $d$,
\item\label{DefAbelianSheafCond4}
$\coker\tau_i$ is annihilated by $\CJ^d$ and $pr_\ast\coker\tau_i$ is a locally free $\CO_S$-module of rank $d$.
\end{enumerate}
A \emph{morphism} between two abelian sheaves $(\CF_i,\Pi_i,\tau_i)$ and $(\CF'_i,\Pi'_i,\tau'_i)$ is a collection of morphisms $\CF_i\to\CF'_i$ which commute with the $\Pi$'s and the $\tau$'s. 
\end{definition}

\noindent
{\it Remark.} 
Abelian sheaves of dimension $d=1$ are called \emph{elliptic sheaves} and were studied by Drinfeld~\cite{Drinfeld3} and Blum-Stuhler~\cite{Blum-Stuhler}. The category of Drinfeld-$A$-modules of rank $r$ over $S$ is anti-equivalent to the category of elliptic sheaves of rank $r$ over $S$ which satisfy $\deg\CF_0=1-r$, see \cite[Theorem 3.2.1]{Blum-Stuhler}.

\begin{proposition}\label{Prop1.6}
The push forward along $\pi:C'_S\to C_S$ defines a \emph{restriction of coefficients functor}
\[
\pi_\ast:\es \UCF'\;=\;(\CF'_i,\Pi'_i,\tau'_i)\es\longmapsto\es\pi_\ast\UCF'\;:=\;(\pi_\ast\CF'_i,\pi_\ast\Pi'_i,\pi_\ast\tau'_i)
\]
from abelian sheaves on $C'$ of rank $r'$, dimension $d'$ and characteristic $c':S\to C'$ over $S$ to abelian sheaves on $C$ of rank $nr'$, dimension $d'$ and characteristic $\pi\circ c':S\to C$ over $S$. Here $n$ is the degree of $\pi$.
\end{proposition}

\begin{proof}
Since $\pi$ is finite and flat the sheaves $\pi_\ast\CF_i$ are locally free of rank $nr'$ by \cite[Corollary 2 to Proposition II.3.2.5]{Bourbaki}. Let $k$ and $\ell$ be relatively prime positive integers with $\frac{k}{\ell}=\frac{d'}{nr'\deg(\infty)}$. Let $e$ be the ramification index of $\pi$ at $\infty'$. Then $n=e\deg(\infty')/\deg(\infty)$ and hence $k=k'/\gcd(k',e)$ and $\ell=\ell'e/\gcd(k',e)$. From axiom~\ref{DefAbelianSheafCond2} of Definition~\ref{Def1.5} we obtain an isomorphism
\[
\Pi'_{i+\ell-1}\circ\ldots\circ\Pi'_i\es:\es\CF'_i\es\isoto\es\CF'_{i+\ell}\otimes_{\CO_{C'_S}}\CO_{C'_S}(-ke\cdot\infty')\,.
\]
Since $\pi^\ast\CO_{C_S}(\infty)=\CO_{C'_S}(e\cdot\infty')$ the projection formula 
\[
\pi_\ast\bigl(\CF'_{i+\ell}\otimes_{\CO_{C'_S}}\CO_{C'_S}(-ke\cdot\infty')\bigr)\es=\es(\pi_\ast\CF'_{i+\ell})\otimes_{\CO_{C_S}}\CO_{C_S}(-k\cdot\infty)
\]
yields
\[
\pi_\ast\Pi'_{i+\ell-1}\circ\ldots\circ\pi_\ast\Pi'_i\es:\es\pi_\ast\CF'_i\es\isoto\es(\pi_\ast\CF'_{i+\ell})\otimes_{\CO_{C_S}}\CO_{C_S}(-k\cdot\infty)
\]
from which the proposition is evident.
\end{proof}

\noindent
{\it Remark.} 
Consider the contravariant moduli functor
\begin{eqnarray*}
\CAB^{r,d}:\es \Sch_{/C}&\longto&\Sets\\[2mm]
(c:S\to C) &\mapsto &\Bigl\{\,\text{Isomorphism classes of abelian sheaves on $C$ of} \\[-2mm]
& & \text{\quad rank $r$, dimension $d$, and characteristic $c$ over $S$}\,\Bigr\}
\end{eqnarray*}
Also this functor is not representable (not even after adding level structures, see \cite[Remark 4.2]{HartlAbSh}). Again the restriction of coefficients functor defines a \emph{restriction of coefficients morphism} 
\[
\pi_\ast:\es\CpAB^{r',d'}\es\longto\es\CAB^{nr',d'}\,,\quad\UCF'\es\mapsto\es\pi_\ast\UCF'\,.
\]

\medskip

\noindent
{\it Remark.} 
If we let $S$ vary, the category of abelian sheaves on $C$ of rank $r$ and dimension $d$ becomes a stack $\AbSh^{r,d}$ for the fppf topology on the category of $\BF_q$-schemes. It is an algebraic stack in the sense of Deligne-Mumford~\cite{DM} by \cite[Theorem 3.1]{HartlAbSh}. The restriction of coefficients functor defines a \emph{restriction of coefficients 1-morphism} $\pi_\ast:\AbShp^{r',d'}\to\AbSh^{nr',d'}$.

The construction of \cite[Theorem 3.2.1]{Blum-Stuhler} yields a 1-isomorphism of $\DrMod^r$ with an open and closed substack of $\AbSh^{r,1}$, see \cite[Example 1.8]{HartlAbSh} such that the following diagram is 2-commutative
\[
\xymatrix @R-0.3pc @M+0.3pc {
\DrModp^{r'} \ar@{^{ (}->}[r] \ar[d]^{\pi_\ast} & \AbShp^{r',1} \ar[d]^{\pi_\ast}\\
\DrMod^{nr'} \ar@{^{ (}->}[r] & \AbSh^{nr',1}\,.
}
\]

\bigskip

Finally let us turn to Drinfeld-Anderson shtuka.

\begin{definition}\label{Def5.1}
A \emph{right} (\emph{left}) \emph{Drinfeld-Anderson shtuka $\ul\CE=(\CE,\wt\CE,j,\tau,b,c)$ on $C$ of rank $r$ and dimension $d$ over $S$} consists of two $\BF_q$-morphisms $b,c:S\to C$ and a diagram
\[
\xymatrix @C=1pc @R=1pc {
\CE \ar[r]^j & \CE' & & & & \CE \\
\sigma^\ast\CE \ar[ur]_\tau& & & \Bigl(\text{resp.} & \CE' \ar[ur]^\tau \ar[r]_j & *!L(0.5) 
\objectbox{\;\sigma^\ast\CE\quad\Bigr)}
}
\]
of locally free sheaves $\CE$ and $\wt\CE$ of rank $r$ on $C_S$ such that $\coker j$, respectively $\coker\tau$, are locally free of rank $d$ as $\CO_S$-modules and supported on the graphs of $b$, respectively $c$. The morphism $b$ is called the \emph{pole of $\ul\CE$} and $c$ is called the \emph{zero of $\ul\CE$}.
\end{definition}

\noindent
{\it Remark.}
Every abelian sheaf $(\CF_i,\Pi_i,\tau_i)$ on $C$ of rank $r$, dimension $d$, and characteristic $c$ over $S$ gives rise to a right Drinfeld-Anderson shtuka on $C$ over $S$ by setting for any $i\in\BZ$
\[
\CE:=\CF_i\,,\quad\wt\CE:=\CF_{i+1}\,,\quad j:=\Pi_i\,,\quad \tau:=\tau_i\,.
\]
This defines a faithful functor from abelian sheaves to Drinfeld-Anderson shtuka on $C$ over $S$. Together with the functor from Drinfeld-$A$-modules to elliptic sheaves on $C$ one obtains a fully faithful functor from Drinfeld-$A$-modules of rank $r$ over $S$ to Drinfeld-Anderson shtuka on $C$ of rank $r$ and dimension $1$ over $S$, see Drinfeld~\cite[\S 1]{Drinfeld5}

\medskip

The argument of Proposition~\ref{Prop1.6} also shows

\begin{proposition}\label{Prop5.2}
Relaxing the conditions on $\pi:C'\to C$ assume only that $\pi$ is finite of degree $n$. Then the push forward along $\pi$ defines a \emph{restriction of coefficients functor}
\[
\pi_\ast:\es (\CE\,,\,\wt\CE\,,\,j\,,\,\tau\,,\,b\,,\,c) \es\longmapsto\es(\pi_\ast\CE\,,\,\pi_\ast\wt\CE\,,\,\pi_\ast j\,,\,\pi_\ast\tau\,,\,\pi\circ b\,,\,\pi\circ c)
\]
from Drinfeld-Anderson shtuka on $C'$ of rank $r'$ and dimension $d'$ to Drinfeld-Anderson shtuka on $C$ of rank $nr'$ and dimension $d'$ over $S$.\qed
\end{proposition}

\noindent
{\it Remark.}
Consider the contravariant moduli functor
\begin{eqnarray*}
\DAS^{r,d}:\es \Sch_{/C\times C}&\longto&\Sets\\[2mm]
\bigl((b,c):S\to C\times_{\BF_q}C) &\mapsto &\Bigl\{\,\text{Isomorphism classes of Drinfeld-Anderson shtuka on} \\[-2mm]
& & \text{\quad $C$ of rank $r$, dimension $d$, pole $b$, and zero $c$ over $S$}\,\Bigr\}
\end{eqnarray*}
Also this functor is not representable but the restriction of coefficients functor defines a \emph{restriction of coefficients morphism} 
\[
\pi_\ast:\es\DASp^{r',d'}\es\longto\es\DAS^{nr',d'}\,,\quad\UCF'\es\mapsto\es\pi_\ast\UCF'\,.
\]
Here again the category of Drinfeld-Anderson shtuka of rank $r$ and dimension $d$ over varying $\BF_q$-schemes $S$ is an algebraic stack $\DrSht^{r,d}$ for the fppf topology in the sense of Deligne-Mumford~\cite{DM} and the restriction of coefficients functor defines a \emph{restriction of coefficients 1-morphism} $\pi_\ast:\DrShtp^{r',d'}\to\DrSht^{nr',d'}$.

 
\section{Restriction of Coefficients for Drinfeld Modules}\label{Sect2}

\begin{theorem}\label{Thm2.1}
The restriction of coefficient morphism $\pi_\ast:\DRMp^{r'}\to\DRM^{nr'}$ for Drinfeld modules is in general not relatively representable.
\end{theorem}

\begin{proof}
We give a counterexample to relative representability. Let $q=3$, $A=\BF_3[x]$, $A'=\BF_3[y]$ and $\pi^\ast\!:A\to A', x\mapsto y^2$. Let $S=\Spec\BF_3$ and $c^\ast\!:A\to\BF_3, x\mapsto 0$. Consider the Drinfeld-$A$-module $(E,\phi)$ of rank $2$ over $S$ given by $E=\BG_{a,S}$ and 
\[
\phi:\es A\longto\BF_3\{\tau\}\,,\quad \phi(x)=\tau^2\,.
\]
Let $\ul T:=\DRMp^1\times_{\DRM^2}S$ be the fiber product of functors. Then $\ul T$ is the contravariant functor
\begin{eqnarray*}
\ul T:\es\Sch_{/\Spec A'\times_{\Spec A}S} & \longto & \Sets \\[2mm]
\Bigl(S',c':S'\to\Spec A'\mbox{ }& \mapsto & \Bigl\{\,\text{Isomorphism classes of Drinfeld-$A'$-modules $(E',\phi')$} \\[-2mm]
f:S'\to S \qquad\mbox{ }\Bigr)& & \text{\quad of rank $1$ over $S'$, such that $f^\ast(E,\phi)\cong\pi_\ast(E',\phi')$}\es\Bigr\}\,.
\end{eqnarray*}
We show that $\ul T$ is not representable. For this purpose make $S$ into a $\Spec A'$-scheme by $(c')^\ast:A'\to\BF_3,y\mapsto0$. Then $\ul T(S)$ contains two isomorphism classes given by $E'_1=E'_2=\BG_{a,S}$ and \[
\phi'_1:\es y\mapsto \tau \qquad\text{and}\qquad \phi'_2:\es y\mapsto-\tau\,.
\]
These two isomorphism classes are different because otherwise there were an isomorphism
\[
\epsilon\es\in\es\Isom\bigl((E'_1,\phi'_1),(E'_2,\phi'_2)\bigr)\es=\es\bigl\{\,\epsilon\in\CO_S\mal:\es-\tau\circ\epsilon=\epsilon\circ\tau\,\bigr\}\,.
\]
That is, $\epsilon\in\BF_3\mal$ must satisfy $-\epsilon^3\tau=\epsilon\tau$, whence $\epsilon^2=-1$. This is impossible for $\epsilon\in\BF_3\mal$.

On the other hand such an element exists in $\BF_9\mal$. So if $S'=\Spec\BF_9$ the two isomorphism classes become equal in $\ul T(S')$. But this implies that $\ul T$ is not representable. Since if it were representable by a scheme $T$ we had two different morphisms from $S$ to $T$ which yield the same morphism from $S'$ to $T$
\[
\Spec\BF_9\es\longto\es\Spec\BF_3\es\raisebox{-.2em}{$\stackrel{\DS\longto}{\DS\longto}$}\es T\,.
\]
As $\Spec\BF_9\to\Spec\BF_3$ is a homeomorphism and $\BF_3\subset\BF_9$ this is impossible.
\end{proof}

\noindent
{\it Remark.} 
The reason why $\ul T$ is not representable is that the isomorphism $\alpha:f^\ast(E,\phi)\isoto\pi_\ast(E',\phi')$ in the definition of $\ul T(S)$ is only supposed to exist but is not added to the data. More precisely we have

\begin{theorem}\label{Thm2.2}
Let $c:S\to\Spec A$ be a morphism of $\BF_q$-schemes and let $(E,\phi)$ be a Drinfeld-$A$-module of rank $nr'$ and characteristic $c$ over $S$. Then the contravariant functor
\begin{eqnarray*}
\ul T:\es\Sch_{/\Spec A'\times_{\Spec A}S} & \longto & \Sets \\[2mm]
\Bigl(S',c':S'\to\Spec A'\mbox{ }& \mapsto & \Bigl\{\,\text{Isomorphism classes of tripples $(E',\phi',\alpha)$ where} \\[-2mm]
f:S'\to S \qquad\mbox{ }\Bigr)& & \quad\bullet\es (E',\phi') \text{ is a Drinfeld-$A'$-module of rank $r'$}\\[-2mm]
& & \qquad\text{ and characteristic $c'$ over $S'$ and}\\[1mm]
& & \quad\bullet\es \alpha:f^\ast(E,\phi)\isoto\pi_\ast(E',\phi') \text{ is a fixed isomorphism}\,\Bigr\}
\end{eqnarray*}
is representable by an affine $S$-scheme of finite presentation.
\end{theorem}

\begin{proof}
Since the question is local on $S$ we may by Proposition~\ref{Prop1.3} assume that $S=\Spec B$, $E=\BG_{a,B}$ and $\phi$ is given by $\phi:A\to B\{\tau\}$ such that the highest coefficient of every $\phi(a)$ is a unit in $B$.

Let the $A$-algebra $A'$ be generated by $a'_1,\ldots,a'_N$. In order to extend $\phi$ to $\phi':A'\to B\{\tau\}$ we must define $\phi'(a'_1),\ldots,\phi'(a'_N)$. Set $d_\nu:=-r'\ord_{\infty'}(a'_\nu)\deg(\infty')$ for all $\nu$. Define
\[
B'\es:=\es B\otimes_A A'\,[\,\delta_{i,\nu},\delta_{d_\nu,\nu}^{-1}:\es\nu=1,\ldots,N\,,\,i=0,\ldots,d_\nu\,]
\]
and the morphism $c':\Spec B'\to\Spec A'$ by the natural map $A'\to B'$. Define
\[
\phi'(a'_\nu)\es:=\es\sum_{i=0}^{d_\nu}\delta_{i,\nu}\tau^i\es\in\es B'\{\tau\}\,,
\]
and $\phi'|_A:=\phi$, and let $\alpha=\id_{\BG_{a,B'}}$. In order that the so defined $\phi'$ is a Drinfeld-$A'$-module of rank $r'$ and characteristic $c'$ over $\Spec B'$ we must require several conditions which are all represented by finitely presented closed subschemes of $\Spec B'$. Namely consider successively for $\nu=1,\ldots,N$ the minimal polynomial of $a'_\nu$ over $A(a'_1,\ldots,a'_{\nu-1})$
\[
(a'_\nu)^m+b_{\nu,m-1}(a'_\nu)^{m-1}+\ldots+b_{\nu,1}a'_\nu+b_{\nu,0}\es=\es0
\]
with $b_{\nu,k}\in A(a'_1,\ldots,a'_{\nu-1})$. The fact that $\phi':A'\to B'\{\tau\}$ is a ring homomorphism is now expressed by the vanishing of
\[
\phi'(a'_\nu)^m+\phi'(b_{\nu,m-1})\phi'(a'_\nu)^{m-1}+\ldots+\phi'(b_{\nu,1})\phi'(a'_\nu)+\phi'(b_{\nu,0})\es=\es0
\]
in $B'\{\tau\}$. Looking at the coefficients of this $\tau$-polynomial we get a finitely generated ideal of $B'$ which we must require to vanish, that is, must divide out. Likewise the commutation of $\phi'(a'_\nu)$ with a (finite) generating system of the $\BF_q$-algebra $A(a'_1,\ldots,a'_{\nu-1})$ yields a finitely generated ideal of $B'$. Finally the condition on the characteristic means that $(c')^\ast(a'_\nu)=\delta_{0,\nu}$. Putting everything together the sum of these ideals defines a closed subscheme $T\subset\Spec B'$ which is of finite presentation and affine over $S$.

We claim that $T$ represents $\ul T$. So let $(E',\phi',\alpha)$ be an element of $\ul T(S')$. The isomorphism $\alpha:f^\ast\BG_{a,S}\isoto E'$ yields an isomorphism $\alpha:(\BG_{a,S'},\psi')\isoto(E',\phi')$ of Drinfeld-$A'$-modules over $S'$ where $\psi'(a):=\alpha^{-1}\circ\phi'(a)\circ\alpha$ for all $a\in A'$. Since $\psi'(a)=f^\ast\phi(a)$ for $a\in A$, $\psi'$ is of the form described in Proposition~\ref{Prop1.3}. In particular
\[
\psi'(a'_\nu)\es=\es\sum_{i=0}^{d_\nu}\delta_i(a'_\nu)\tau^i\es\in\es\Gamma(S',\CO_{S'})\{\tau\}\,.
\]
Mapping $\delta_{i,\nu}$ to $\delta_i(a'_\nu)$ defines the desired uniquely determined morphism $B'\to\Gamma(S',\CO_{S'})$, whence $S'\to T$.
\end{proof}

\begin{proposition}\label{Prop2.3}
In the situation of Theorem~\ref{Thm2.2} the scheme $T$ representing $\ul T$ is finite over $S$.
\end{proposition}

\begin{proof}
We already know that $T$ is separated and of finite presentation over $S$. We use the valuative criterion of properness to show that it is proper. Since it is also affine over $S$ it must be finite.

So let $R$ be a valuation ring with fraction field $K$ and consider the diagram
\[
\xymatrix @R=2.5pc @C=6pc{
\Spec K \ar[r]^{(E',\phi',\alpha)} \ar[d]_g & T \ar[d]^{\pi_\ast} \\
\Spec R \ar[r]_f \ar@{-->}[ru]_{(\wt E,\wt \phi,\wt\alpha)} & S
}
\]
where the horizontal arrow on top is given by a Drinfeld-$A'$-module $(E',\phi')$ of rank $r'$ and characteristic $c':\Spec K\to\Spec A'$ 
together with an isomorphism \mbox{$\alpha:(fg)^\ast(E,\phi)\isoto\pi_\ast(E',\phi')$} over $\Spec K$. We must exhibit the dashed arrow which corresponds to a Drinfeld-$A'$-module $(\wt E,\wt\phi)$ of rank $r'$ and characteristic $\tilde c:\Spec R\to\Spec A'$ (note that $c'$ factors through a unique morphism $\tilde c$ satisfying $\pi\circ\tilde c=c\circ f:\Spec R\to \Spec A$ because $\Spec A'$ is proper over $\Spec A$) together with an isomorphism $\wt \alpha:f^\ast(E,\phi)\isoto\pi_\ast(\wt E,\wt\phi)$ over $\Spec R$. The commutativity of the diagram means that there exists an isomorphism $\beta:g^\ast(\wt E,\wt\phi)\isoto(E',\phi')$ with $\pi_\ast\beta\circ g^\ast\wt\alpha=\alpha$.

Since $R$ is a local ring $f^\ast E=\BG_{a,R}$ without loss of generality and $f^\ast\phi:A\to R\{\tau\}$. We use the isomorphism $\alpha$ to replace $(E',\phi')$ by $(\BG_{a,K},\psi')$ with $\psi'(a):=\alpha^{-1}\circ\phi'(a)\circ\alpha\in K\{\tau\}$ for all $a\in A'$. Thus $\alpha$ is replaced by $\id_{\BG_{a,K}}$ and $\psi'(a)=(fg)^\ast\phi(a)$ for all $a\in A$. If we show that $\psi'(a)$ belongs to $R\{\tau\}$ for all $a\in A'$, then we may take $\wt E=\BG_{a,R}$ and $\wt\phi=\psi':A'\to R\{\tau\}$, as well as $\wt\alpha=\id_{\BG_{a,R}}$ and $\beta=\id_{\BG_{a,K}}$, and we are done.

So let $a\in A'$ and let $a^m+b_{m-1}a^{m-1}+\ldots+b_1a+b_0=0$ be an equation of integral dependence of $a$ over $A$. In particular
\begin{equation}\label{Eq2.1}
\psi'(a)^m+\psi'(b_{m-1})\psi'(a)^{m-1}+\ldots+\psi'(b_1)\psi'(a)+\psi'(b_0)\es=\es0\,.
\end{equation}
Over an algebraic closure of $K$ the polynomial $\psi'(a)(x)$, where we use $\tau(x)=x^q$, splits as $\psi'(a)(x)=\prod_i (x-\lambda_i)$ with $\lambda_i\in K^\alg$. From equation (\ref{Eq2.1}) we see that each $\lambda_i$ is a root of $\psi'(b_0)=(fg)^\ast\phi(b_0)$. Since $f^\ast\phi(b_0)$ has coefficients in $R$ with the highest coefficient in $R\mal$, all $\lambda_i$ must be integral over $R$. Therefore the coefficients of $\psi'(a)$ which are symmetric polynomials in the $\lambda_i$ are integral over $R$ and belong to $K$, hence they lie in $R$ as desired. This proves the proposition.
\end{proof}

\begin{theorem}\label{Thm2.4}
The restriction of coefficients morphism $\pi_\ast:\DRMp^{r'}\to\DRM^{nr'}$ satisfies the valuative criterion for properness.
\end{theorem}

\begin{proof}
Let $R$ be a valuation ring with fraction field $K$ and consider the diagram
\[
\xymatrix @C+2pc @R+1pc @M+0.2pc {
 & \ul T \ar[d]\\
\Spec K\ar@/^/[ru]^{(E',\phi',\alpha)} \ar[r]^{(E',\phi')\quad\qquad\qquad\mbox{ }} \ar[d]_f & 
**{=<15pc,2pc>} \objectbox{\DRMp^{r'}\times_{\DRM^{nr'}}\Spec R} \ar[d]^{\qquad\qquad\qquad\quad\mbox{\large$\Box$}} \ar[r] & \DRMp^{r'} \ar[d]^{\pi_\ast}\\
\Spec R \ar@{-->}[ru] \ar@{=}[r] & \Spec R \ar[r]^{(E,\phi)} & \DRM^{nr'}
}
\]
where the horizontal morphisms are induced by a Drinfeld-$A'$-module $(E',\phi')$ of rank $r'$ and characteristic $c':\Spec K\to\Spec A'$ over $\Spec K$ and a Drinfeld-$A$-module $(E,\phi)$ of rank $nr'$ and characteristic $c:\Spec R\to\Spec A$ over $\Spec R$ and where $\ul T$ is the representable functor from Theorem~\ref{Thm2.2} for $S=\Spec R$. The commutativity of the square on the left means that $f^\ast(E,\phi)\cong\pi_\ast(E',\phi')$. The choice of any such isomorphism $\alpha$ defines a morphism $\Spec K\to \ul T$. By Proposition~\ref{Prop2.3} we find a unique morphism $\Spec R\to \ul T$ fitting into the diagram which induces the dashed morphism.
It remains to show that the dashed morphism is uniquely determined (independent of the choice of $\alpha$) and this is proved in the following lemma.
\end{proof}

\begin{lemma}\label{Lemma2.5}
Let $R$ be a valuation ring with fraction field $K$ and let $f:\Spec K\to\Spec R$ be the induced morphism. Let $(E'_1,\phi'_1)$ and $(E'_2,\phi'_2)$ be two Drinfeld-$A'$-modules of rank $r'$ and characteristic $c':\Spec R\to\Spec A'$ over $\Spec R$ and let $\alpha:f^\ast(E'_1,\phi'_1)\isoto f^\ast(E'_2,\phi'_2)$ be an isomorphism over $\Spec K$. Then $\alpha=f^\ast\beta$ for a unique isomorphism $\beta:(E'_1,\phi'_1)\isoto (E'_2,\phi'_2)$ over $\Spec R$.
\end{lemma}

\begin{proof}
Since $R$ is a local ring we have without loss of generality $E'_1=E'_2=\BG_{a,R}$. Let $a\in A'\setminus\BF_q$ and write for $j=1,2$
\[
\phi'_j(a)\es=\es\sum_{i=0}^m \delta_{i,j}\tau^i\,.
\]
The isomorphism over $\Spec K$ is given by an element $\alpha\in K\mal$ which satisfies $\phi'_2(a)\circ\alpha=\alpha\circ\phi'_1(a)$, whence $\delta_{2,m}\,\alpha^{q^m}=\alpha\,\delta_{1,m}$. Since $\delta_{1,m}$ and $\delta_{2,m}$ are units in $R$ the same is true for $\alpha$. So the isomorphism $\alpha$ is already defined over $\Spec R$. 
\end{proof}

\noindent
{\it Remark.} 
Phrased in the language of stacks \cite{LaumonMB}, Theorems~\ref{Thm2.1} and \ref{Thm2.4} say that the restriction of coefficients 1-morphism $\pi_\ast:\DrModp^{r'}\to\DrMod^{nr'}$ is proper but in general not representable. Namely by the arguments of Theorem~\ref{Thm2.2} the stack $\CT$ classifying data $\bigl((E,\phi),\,(E',\phi,),\,\alpha\bigr)$ where $(E,\phi)$, respectively $(E',\phi')$, is a Drinfeld-$A$-module of rank $nr'$, respectively a Drinfeld-$A'$-module of rank $r'$ over the same scheme $S$ together with a fixed isomorphism $\alpha:(E,\phi)\isoto\pi_\ast(E',\phi')$ over $S$  is relatively representable over $\DrMod^{nr'}$ by a finite and finitely presented morphism of schemes. The projection $\CT\to\DrModp^{r'}$ onto $(E',\phi')$ is an \'etale epimorphism and makes $\CT$ into a torsor under the finite relative group scheme $\Aut\bigl(\pi_\ast(E',\phi')\bigr)$ over $\DrModp^{r'}$. In particular $\DrModp^{r'}$ is of finite presentation over $\DrMod^{nr'}$ since $\CT$ is and it satisfies the valuative criterion for properness by the arguments of Theorem~\ref{Thm2.4}.


\section{Restriction of Coefficients for Abelian Sheaves} \label{Sect3}

\begin{theorem}\label{Thm3.1}
The restriction of coefficients morphism $\pi_\ast:\CpAB^{r',d'}\to\CAB^{nr',d'}$ for abelian sheaves is in general not relatively representable.
\end{theorem}

\begin{proof}
This follows directly from Theorem~\ref{Thm2.1} and the remark after Definition~\ref{Def1.5}. The example from Theorem~\ref{Thm2.1} yields the following abelian sheaf on $C=\BP^1_{\BF_3}$ over $S=\Spec\BF_3$. Let $\CF_i=\CO_{C_S}\bigl(\lfloor\frac{i+1}{2}\rfloor\cdot\infty\bigr)\oplus\CO_{C_S}\bigl(\lfloor\frac{i}{2}\rfloor\cdot\infty\bigr)$ where $\lfloor\frac{i}{2}\rfloor$ is the largest integer $\le\frac{i}{2}$. Let $\Pi_i$ be the natural inclusion $\CF_i\subset\CF_{i+1}$ and let $\tau_i:\sigma^\ast\CF_i\to\CF_{i+1}$ be given by the matrix $\left(\begin{array}{cc}0 & x \\ 1 & 0 \end{array}\right)$ where $\BP^1_{\BF_3}\setminus\{\infty\}=\Spec\BF_3[x]$.

Let $\pi:C'=\BP^1_{\BF_3}\to C$ be given by $A\to A'=\BF_3[y]$, $x\mapsto  y^2$. Then $(\CF_i,\Pi_i,\tau_i)$ is isomorphic to $\pi_\ast(\CF'_i,\Pi'_i,\tau'_i)$ for $\CF'_i=\CO_{C'_S}(i\cdot\infty')$, $\Pi'_i$ the natural inclusion, and $\tau'_i=\pm y$. The two abelian sheaves for $\tau'_i=+y$ and $\tau'_i=-y$ are not isomorphic over $\Spec\BF_3$ but become isomorphic over $\Spec\BF_9$.
\end{proof}

\begin{theorem}\label{Thm3.2}
Let $S$ be a locally noetherian $\BF_q$-scheme and let $c:S\to C$ be an $\BF_q$-morphism. Let $\UCF$ be an abelian sheaf on $C$ of rank $nr'$, dimension $d'$ and characteristic $c$ over $S$. Then the contravariant functor
\begin{eqnarray*}
\ul T:\es\Sch_{/C'\times_{C}S} & \longto & \Sets \\[2mm]
\Bigl(S',c':S'\to C'\mbox{ }& \mapsto & \Bigl\{\,\text{Isomorphism classes of pairs $(\UCF',\alpha)$ where} \\[-2mm]
\mbox{ }\qquad f:S'\to S \mbox{ }\Bigr)& & \quad\bullet\es \UCF' \text{ is an abelian sheaf of rank $r'$, dimension $d'$,}\\[-2mm]
& & \qquad\text{ and characteristic $c'$ over $S'$ and}\\[1mm]
& & \quad\bullet\es \alpha:f^\ast\UCF\isoto\pi_\ast\UCF' \text{ is a fixed isomorphism}\,\Bigr\}
\end{eqnarray*}
is representable by a (quasi-affine) $S$-scheme of finite type.
\end{theorem}

For the proof we need the following 

\begin{lemma}\label{Lemma3.3}
Let $S$ be a locally noetherian scheme, let $\rho:Y\to S$ be a flat projective morphism, and let $\pi:X\to Y$ be a finite faithfully flat morphism of degree $n$. For an $S$-scheme $S'$ set $Y':=Y\times_S S'$ and $X':=X\times_S S'$. Let $\CF$ be a locally free sheaf on $Y$ of rank $rn$. Then the contravariant functor
\begin{eqnarray*}
\ul U:\es\Sch_{/S} & \longto & \Sets \\[2mm]
(f:S'\to S )& \mapsto & \Bigl\{\,\text{Isomorphism classes of pairs $(\CF',\alpha)$ where} \\[-1mm]
& & \quad\bullet\es \CF' \text{ is a locally free sheaf of rank $r$ on $X'$ and}\\[0mm]
& & \quad\bullet\es \alpha:f^\ast\CF\isoto\pi_\ast\CF' \text{ is a fixed isomorphism}\,\Bigr\}
\end{eqnarray*}
is representable by a (quasi-affine) $S$-scheme of finite type.
\end{lemma}

\begin{proof}
Since the question is local on $S$ we may assume that $S$ is affine. By \cite[II, Proposition 1.4.3]{EGA} the functor $\ul U$ is isomorphic to the functor
\[
\ul U':\es (f:S'\to S)\es\mapsto\es\Hom_{\CO_{Y'}\text{-algebras}}\bigl(\pi_\ast\CO_{X'}\,,\,\CE nd_{\CO_{Y'}} (f^\ast\CF)\bigr)
\]
the set of $\CO_{Y'}$-algebra homomorphisms $\pi_\ast\CO_{X'}\to\CE nd_{\CO_{Y'}} (f^\ast\CF)$. Fix an ample invertible sheaf $\CL$ on $Y$. For any integer $N$ define $\CH_N\;:=\;\CHom_{\CO_Y}\bigl(\CF,\CF\otimes_{\CO_Y}\CL^{\otimes N}\bigr)\;=\;\CE nd_{\CO_Y} (f^\ast\CF)\otimes_{\CO_Y}\CL^{\otimes N}$. Then for the homomorphisms of $\CO_{Y'}$-\emph{modules} we obtain
\[
\Hom_{\CO_{Y'}\text{-modules}}\bigl(\pi_\ast\CO_{X'}\,,\,\CE nd_{\CO_{Y'}}(f^\ast\CF)\bigr)\es=\es\Hom_{\CO_{Y'}\text{-mod}}\bigl(\pi_\ast\CO_{X'}\otimes_{\CO_{Y'}}\CL^{\otimes N}\,,\,f^\ast\CH_N\bigr)\,.
\]
There is an integer $N$ such that
\begin{itemize}
\item 
$\pi_\ast\CO_X\otimes_{\CO_Y}\CL^{\otimes N}$ is generated by global sections and
\item 
$\rho_\ast(\pi_\ast\CO_X\otimes_{\CO_Y}\CL^{\otimes N})$ and $\rho_\ast\CH_N$ are locally free on $S$
\end{itemize}
since these conditions are achieved for $N\gg0$ by the Theorem on Cohomology and Base Change \cite[Theorem III.12.11]{Hartshorne}.

Shrinking $S$ we let $x_1,\ldots,x_m$ be an $\CO_S$-basis of $\rho_\ast(\pi_\ast\CO_X\otimes_{\CO_Y}\CL^{\otimes N})$. We must specify their images in $\rho_\ast\CH_N$. Let $U_1:=\ul\Spec_S\Sym_{\CO_S}(\rho_\ast\CH_N)\dual$ and $U_2:=U_1\times_S\ldots\times_S U_1$ the $m$-fold fiber product. Let $f:U_2\to S$ be the induced morphism and set $Y_2:=Y\times_S U_2$ and $X_2:=X\times_S U_2$. Then for any $S$-scheme $S'$
\begin{eqnarray*}
\Hom_S(S',U_1) & = & \Hom_{\CO_S\text{-algebras}}\bigl( \Sym_{\CO_S}(\rho_\ast\CH_N)\dual\,,\,\CO_{S'}\bigr) \\[2mm]
& = & \Hom_{\CO_S\text{-modules}}\bigl((\rho_\ast\CH_N)\dual\,,\,\CO_{S'}\bigr) \\[2mm]
& = & \Gamma\bigl(S'\,,\,\CO_{S'}\otimes_{\CO_S}\rho_\ast\CH_N\bigr)\,.
\end{eqnarray*}
So on $U_2$ there exist $m$ universal global sections of $\rho_\ast\CH_N$ which we use as the images of our $x_1,\ldots,x_m$ to obtain a universal homomorphism of $\CO_{U_2}$-modules
\begin{equation}\label{Eq3.1}
\rho_\ast\bigl(\pi_\ast\CO_{X_2}\otimes_{\CO_{Y_2}}f^\ast\CL^{\otimes N}\bigr)\es\longto\es\rho_\ast f^\ast\CH_N\,.
\end{equation}

Next we take care of the $\CO_Y$-algebra structures. Every $x_i$ has a minimal polynomial over $\Gamma(Y,\CL^{\otimes N})[x_1,\ldots,x_{i-1}]$ of the form
\[
P(x_i)\es:=\es x_i^k+a_{k-1}x_i^{k-1}+\ldots+a_0\es=\es0
\]
inside $\Gamma(Y,\pi_\ast\CO_X\otimes_{\CO_Y}\CL^{\otimes Nk})$. Using our homomorphism (\ref{Eq3.1}) and the $\CO_{U_2}$-module structure of $f^\ast\CF$ we view $P(x_i)$ as an element of $\Gamma(U_2,f^\ast\rho_\ast\CH_{Nk})$. The requirement that this element vanishes defines a closed subscheme of $U_2$ by Lemma~\ref{Lemma3.4} below. Let $U_3$ be the closed subscheme of $U_2$ obtained in this way for $i=1,\ldots,m$. This yields a $\pi_\ast\CO_{X_3}$-module structure on $f^\ast\CF$, whence (an isomorphism class of) a coherent sheaf $\CF_3$ on $X_3:=X\otimes_S U_3$ together with an isomorphism $\alpha:f^\ast\CF\isoto \pi_\ast\CF_3$.

It remains to represent the condition that $\CF_3$ is locally free. Let $V\subset X_3$ be the open subscheme on which $\CF_3$ is flat, see \cite[IV$_3$, Theorem 11.1.1]{EGA}. Define $U:=U_3\setminus \pi(X_3\setminus V)$. Since $\rho\pi:X_3\to U_3$ is proper $U\subset U_3$ is open. Since $(\rho\pi)^{-1}U\subset V$ the coherent sheaf $\CF_3$ is locally free on $(\rho\pi)^{-1}U$ of rank $r$. We claim that $U$ represents the functor $\ul U$. Indeed, let $S'$ be an $S$-scheme and $(\CF',\alpha)\in \ul U(S')$. Then the $\pi_\ast\CO_{X'}$-module structure on $\pi_\ast\CF'$ defines a uniquely determined morphism $S'\to U_3$. Since above every point $s\in S'$ the fiber $\CF'_s$ is flat on $X\times_S s$, the image of $s$ in $U_3$ lands in $U$ by \cite[IV$_3$, Theorem 11.3.10]{EGA}. (This is the only place where we use the assumption that $\pi$ is flat.) This proves the lemma.
\end{proof}

\begin{lemma}\label{Lemma3.4}
Let $S$ be a scheme and let $\CH$ be a locally free sheaf on $S$. Let $I$ be a set and let $h_i\in \Gamma(S,\CH)$ for all $i\in I$. Then the condition $h_i=0$ for all $i\in I$ is represented by a closed subscheme of $S$.
\end{lemma}

\begin{proof}
This is \cite[0$_{\rm new}$, Proposition 5.5.1]{EGA} taking into account that on a locally noetherian topological space the set of global sections of an arbitrary direct sum equals the direct sum of the global sections.
\end{proof}

\begin{proof}[Proof of Theorem~\ref{Thm3.2}]
Let $\UCF=(\CF_i,\Pi_i,\tau_i)$ and let $\ell'$ and $k'$ be relatively prime positive integers with $\frac{k'}{\ell'}=\frac{d'}{r'\deg(\infty')}$. For $i=0,\ldots,\ell'$ let $U_i$ be the scheme from Lemma~\ref{Lemma3.3} classifying the pairs $(\CF'_i,\alpha_i)$ of locally free sheaves $\CF'_i$ on $X=C'_S$ and isomorphisms $\alpha_i:\CF_i\isoto\pi_\ast\CF'_i$. Over $T:=U_0\times_S\ldots\times_S U_{\ell'}$ we have the universal sheaves $\CF'_0,\ldots,\CF'_{\ell'}$ on $C'_T$. We need that the morphisms of $\CO_{C_T}$-modules
\begin{eqnarray*}
\Pi'_i\es:=&\es\:\alpha_{i+1}\circ\Pi_i\circ\alpha_i^{-1}:&\pi_\ast\CF'_i\es\longto\es\pi_\ast\CF'_{i+1} \qquad\text{and}\\[2mm] 
\tau'_i\es:=&\alpha_{i+1}\circ\tau_i\circ\sigma^\ast\alpha_i^{-1}:&\pi_\ast\sigma^\ast\CF'_i\es\longto\es\pi_\ast\CF'_{i+1}
\end{eqnarray*}
are actually morphisms of $\pi_\ast\CO_{C'_T}$-modules and thus by \cite[II, Proposition 1.4.3]{EGA} morphisms $\Pi':\CF'_i\to\CF'_{i+1}$ and $\tau'_i:\sigma^\ast\CF'_i\to\CF'_{i+1}$.

It suffices to work on an affine covering of $T$. Let $pr:C_T\to T$ be the projection onto the second factor. Let $\CL$ be an ample invertible sheaf on $C$ and let $N$ be an integer such that for $i=0,\ldots, \ell'-1$
\begin{itemize}
\item 
$\pi_\ast\CO_{C'_T}\otimes_{\CO_{C_T}}\CL^{\otimes N}$ is generated by global sections $x_1,\ldots,x_m$, 
\item 
$\pi_\ast\CF'_i\otimes_{\CO_{C_T}}\CL^{\otimes N}$ is generated by global sections $y_1,\ldots,y_n$, and
\item 
$\CH_{i+1}:=pr_\ast(\pi_\ast\CF'_{i+1}\otimes_{\CO_{C_T}}\CL^{\otimes 2N})$ is locally free on $T$.
\end{itemize}
Then $\CG_i:=\pi_\ast\CO_{C'_T}\otimes_{\CO_{C_T}}\pi_\ast\CF'_i\otimes_{\CO_{C_T}}\CL^{\otimes 2N}$ is generated by the $x_\mu\otimes y_\nu$. There are two morphisms of $\CO_{C_T}$-modules
\[
\CG_i \es\raisebox{-.2em}{$\stackrel{\DS\longto}{\DS\longto}$}\es\pi_\ast\CF'_{i+1}\otimes_{\CO_{C_T}}\CL^{\otimes 2N}
\]
depending on the order in which $\Pi'_i$ is composed with the contraction $\pi_\ast\CO_{C'_T}\otimes_{\CO_{C_T}}\pi_\ast\CF'_i\to\pi_\ast\CF'_i$ (coming from the $\CO_{C'_T}$-module structure on $\CF'_i$). Whether the difference of these two morphisms is the zero morphism can be tested on the images of the global sections $x_\mu\otimes y_\nu$ inside $\CH_{i+1}$. By Lemma~\ref{Lemma3.4} this condition is represented by a closed subscheme of $T$.

We proceed analogously for the $\tau_i$ and obtain a closed subscheme $T_1\subset T$ and for $i=0,\ldots,\ell'-1$ universal morphisms $\Pi'_i:\CF'_i\to\CF'_{i+1}$ and $\tau'_i:\sigma^\ast\CF'_i\to\CF'_{i+1}$ on $C'_{T_1}$ which satisfy axiom~\ref{DefAbelianSheafCond1} of Definition~\ref{Def1.5}. Since $pr_\ast\pi_\ast\coker\Pi'_i=pr_\ast\coker\Pi_i$, and the same for $\tau_i$, also axioms~\ref{DefAbelianSheafCond3} and \ref{DefAbelianSheafCond4} hold except for the condition on the support. For this condition let $T_2:=C'\times_C T_1$, let $c':T_2\to C'$ be the projection and let $\CJ'$ be the ideal defining the graph of $c'$. Similarly to the above argument let $\CL$ and $N$ be such that $(\CJ')^{\otimes d'}\otimes_{\CO_{C'_{T_2}}}\CL^{\otimes N}$ is generated by global sections. Again by Lemma~\ref{Lemma3.4} the condition that the multiplication morphism
\[
pr_\ast\bigl((\CJ')^{\otimes d'}\otimes_{\CO_{C'_{T_2}}}\CL^{\otimes N}\otimes_{\CO_{C'_{T_2}}}\coker\tau'_i\bigr)\es\longto\es pr_\ast\bigl(\CL^{\otimes N}\otimes_{\CO_{C'_{T_2}}}\coker\tau'_i\bigr)
\]
is zero is represented by a closed subscheme $T_3$ of $T_2$.

Finally for axiom~\ref{DefAbelianSheafCond2} consider the morphism
\begin{equation} \label{Eq3.2}
\Pi'_{\ell'-1}\circ\ldots\circ\Pi'_0:\es\CF'_0\es\longto\es\CF'_{\ell'}\es\longto\es\CF'_{\ell'}\otimes_{\CO_{C'_{T_3}}}\,\CO_{C'_{T_3}}/\CO_{C'_{T_3}}(-k'\cdot\infty')\,.
\end{equation}
Since $\coker\Pi'_i$ has rank $d'$ axiom~\ref{DefAbelianSheafCond2} is satisfied if and only if the morphism (\ref{Eq3.2}) is the zero morphism. Using that the target is locally free on $T_3$ and reasoning as above the later condition is represented by a closed subscheme $T_4$ of $T_3$. Over $T_4$ we define $\CF'_{i+m\ell'}:=\CF'_i(k'm\cdot\infty')$ for all $i=0,\ldots,\ell'-1$ and all $m\in \BZ$. Then $T_4$ represents the functor $\ul T$.
\end{proof}

\begin{proposition}\label{Prop3.5}
In the situation of Theorem~\ref{Thm3.2} the scheme $T$ representing $\ul T$ is finite over $S$.
\end{proposition}

\begin{proof}
By Theorem~\ref{Thm3.2} it is separated, of finite type, and quasi-affine over $S$. It remains to show that $T$ is proper over $S$. So let $R$ be a discrete valuation ring with fraction field $K$ and consider the diagram
\[
\xymatrix @R=2.5pc @C=6pc{
\Spec K \ar[r]^{(\UCF',\alpha)} \ar[d]_g & T \ar[d]^{\pi_\ast} \\
\Spec R \ar[r]_f \ar@{-->}[ru]_{(\ul{\wt\CF},\wt\alpha)} & S
}
\]
where the horizontal arrow is given by an abelian sheaf $\UCF'$ on $C'$ over $\Spec K$ of rank $r'$, dimension $d'$ and characteristic $c':\Spec K\to C'$ together with an isomorphism $\alpha:(fg)^\ast\UCF\isoto\pi_\ast\UCF'$ on $C_K$. We need to construct an abelian sheaf $\ul{\wt\CF}$ on $C'$ over $\Spec R$ of rank $r'$, dimension $d'$, and characteristic $\tilde c:\Spec R\to C'$ (again the properness of $\pi$ implies that $c'$ factors through a unique morphism $\tilde c$ with $\pi\circ\tilde c=c\circ f:\Spec R\to C$) together with an isomorphism $\wt\alpha:f^\ast\UCF\isoto\pi_\ast\ul{\wt\CF}$ on $C_R$ and an isomorphism $\beta:g^\ast\ul{\wt\CF}\isoto\UCF'$ on $C'_K$ satisfying $\pi_\ast\beta\circ g^\ast\wt\alpha=\alpha$.

We begin by constructing for all $i\in\BZ$ the locally free sheaf $\wt\CF_i$ on $C'_R$ and the isomorphism $\wt\alpha_i:f^\ast\CF_i\isoto\pi_\ast\wt\CF_i$. Let $\CL$ be an ample invertible sheaf on $C$ such that $\pi_\ast\CO_{C'}\otimes_{\CO_C}\CL$ is generated by global sections $x_1,\ldots,x_m$. 

For the next step in the proof we need to introduce some notation. Let $\varpi$ be the generic point of the special fiber of $C_R$ over the residue field of $R$ and let $\CO_\varpi:=\CO_{C_R,\varpi}$ be the local ring at $\varpi$. It is a discrete valuation ring and every uniformizing parameter of $R$ is a uniformizing parameter of $\CO_\varpi$. Let further $K(C)$ be the fraction field of $\CO_\varpi$. It equals the function field of $C_K$. Similarly let $\CO_{\varpi'}$ and $K(C')$ be the rings associated with the curve $C'$. Since the $f^\ast\Pi_i$ are invertible over $\CO_\varpi$ we get $\tau$-modules $\bigl(f^\ast\CF_i\otimes_{\CO_{C_R}}\CO_\varpi,f^\ast(\Pi_i^{-1}\circ\tau_i)\bigr)$ over $\CO_\varpi$ with $f^\ast(\Pi_i^{-1}\circ\tau_i)$ being isomorphisms. Now the argument of Gardeyn~\cite[Proposition 2.13({\it i}\,)]{Gardeyn4} shows that $\bigl(f^\ast\CF_i\otimes_{\CO_{C_R}}\CO_\varpi,f^\ast(\Pi_i^{-1}\circ\tau_i)\bigr)$ is the unique maximal $f^\ast(\Pi_i^{-1}\circ\tau_i)$-invariant $\CO_\varpi$-lattice in $\bigl(f^\ast\CF_i\otimes_{\CO_{C_R}}K(C),f^\ast(\Pi_i^{-1}\circ\tau_i)\bigr)$. Since every $x_\mu$ is an endomorphism of  
\[
\bigl(f^\ast\CF_i\otimes_{\CO_{C_R}}K(C)\,,\,f^\ast(\Pi_i^{-1}\circ\tau_i)\bigr)
\]
it must map $f^\ast\CF_i\otimes_{\CO_{C_R}}\CO_\varpi$ into itself. This makes $f^\ast\CF_i\otimes_{\CO_{C_R}}\CO_\varpi$ into a free $\CO_{\varpi'}$-module.
Now we can apply Lafforgue's \cite[Lemme 2.7]{Lafforgue1} which says that to give a locally free sheaf $\wt\CF_i$ on $C'_R$ is equivalent to giving its restrictions $\wt\CF_i\otimes_{\CO_{C'_R}}\CO_{C'_K}$ and $\wt\CF_i\otimes_{\CO_{C'_R}}\CO_{\varpi'}$. Thus out of $\CF'_i$ and the $\CO_{\varpi'}$-module $f^\ast\CF_i\otimes_{\CO_{C_R}}\CO_\varpi$ we may construct the locally free sheaf $\wt\CF_i$ together with the isomorphism $\beta_i:g^\ast\wt\CF_i\isoto\CF'_i$. Since by construction
\[
\alpha_i:\es\Bigl((fg)^\ast\CF_i\,,\,f^\ast\CF_i\otimes_{\CO_{C_R}}\CO_\varpi\Bigr)\es\isoto\es\Bigl(\pi_\ast(\wt\CF_i\otimes_{\CO_{C'_R}}\CO_{C'_K})\,,\,\pi_\ast(\wt\CF_i\otimes_{\CO_{C'_R}}\CO_{\varpi'})\Bigr)
\]
is an isomorphism on the two restrictions we obtain the isomorphism $\wt\alpha_i:f^\ast\CF_i\isoto\pi_\ast\wt\CF_i$ from Lafforgue's lemma.

Since the $\Pi'_i$ and the $\tau'_i$ are commuting homomorphisms of $\CO_{C'_K}$-modules they restrict to commuting homomorphisms $\wt\Pi_i$ and $\wt\tau_i$ of $\CO_{C'_R}$-modules. Altogether we have shown that $\ul{\wt\CF}=(\wt\CF_i,\wt\Pi_i,\wt\tau_i)$ satisfies axioms~\ref{DefAbelianSheafCond1}, \ref{DefAbelianSheafCond3}, and \ref{DefAbelianSheafCond4} from Definition~\ref{Def1.5} except for the condition on the support of $\coker\wt\tau_i$. Let $\wt\CJ$ be the ideal sheaf on $C'_R$ defining the graph of $\tilde c$. Then $\wt\CJ^{d'}$ annihilates the generic fiber of the free $R$-module $\coker\wt\tau_i$, so it annihilates all of $\coker\wt\tau_i$. Likewise if $z'$ is a uniformizing parameter on $C'$ at $\infty'$ then $(z')^{k'}$ annihilates the generic fiber of the free $R$-module $\coker(\wt\Pi_{i+\ell'-1}\circ\ldots\circ\wt\Pi_i)$, so it annihilates this whole cokernel. Now all axioms are verified and $\ul{\wt\CF}$ is the desired abelian sheaf on $C'$ over $\Spec R$.
\end{proof}

\begin{theorem}\label{Thm3.6}
The restriction of coefficients morphism $\pi_\ast:\CpAB^{r',d'}\to\CAB^{nr',d'}$ satisfies the valuative criterion for properness.
\end{theorem}

\begin{proof}
Since the stacks $\AbSh^{r,d}$ are locally noetherian by \cite[Theorem 3.1]{HartlAbSh} it suffices to test the valuative criterion only for \emph{discrete} valuation rings. For those the argument proceeds as in Theorem~\ref{Thm2.4} using Lemma~\ref{Lemma3.7} below instead of Lemma~\ref{Lemma2.5}.
\end{proof}

\begin{lemma}\label{Lemma3.7}
Let $R$ be a valuation ring with fraction field $K$ and let $f:\Spec K\to\Spec R$ be the induced morphism. Let $\ul\CF$ and $\ul\CF'$ be two abelian sheaves on $C$ over $\Spec R$ of rank $r$, dimension $d$, and characteristic $c:\Spec R\to C$. Let $\alpha:f^\ast\ul\CF\to f^\ast\ul\CF'$ be an isomorphism over $\Spec K$. Then $\alpha=f^\ast\beta$ for a unique isomorphism $\beta:\ul\CF\isoto\ul\CF'$ over $\Spec R$.
\end{lemma}

\begin{proof}
Recall the rings $\CO_{\varpi}$ and $K(C)$ introduced in the proof of Proposition~\ref{Prop3.5} and consider the $\tau$-modules $\bigl(\CF_i\otimes_{\CO_{C_R}}\CO_\varpi,\Pi_i^{-1}\circ\tau_i\bigr)$ and $\bigl(\CF'_i\otimes_{\CO_{C_R}}\CO_\varpi,\Pi_i'{}^{-1}\circ\tau'_i\bigr)$ over $\CO_\varpi$. By the arguments of Gardeyn~\cite[Proposition 2.13({\it i}\,)]{Gardeyn4} these are the unique maximal $\Pi_i^{-1}\circ\tau_i$-invariant $\CO_\varpi$-modules in \mbox{$\CF_i\otimes_{\CO_{C_R}}K(C)$}, respectively $\CF'_i\otimes_{\CO_{C_R}}K(C)$. Hence they are mapped isomorphically into each other under the isomorphism $\alpha$. Now the lemma follows from \cite[Lemme 2.7]{Lafforgue1}.
\end{proof}

\noindent
{\it Remark.}
Like for Drinfeld modules these results say in the language of stacks that the restriction of coefficients 1-morphism $\pi_\ast:\AbShp^{r',d'}\to\AbSh^{nr',d'}$ is proper but in general not representable.


\section{Restriction of Coefficients for Drinfeld-Anderson Shtuka} \label{Sect4}

\begin{theorem}\label{Thm5.3}
The restriction of coefficients morphism $\pi_\ast:\DASp^{r',d'}\to\DAS^{nr',d'}$ for Drinfeld-Anderson shtuka is in general not relatively representable.
\end{theorem}

\begin{proof}
The abelian sheaf from Theorem~\ref{Thm3.1} provides the counter example also for Drinfeld-Anderson shtuka.
\end{proof}

The same reasoning as in Theorem~\ref{Thm3.2} and Proposition~\ref{Prop3.5} yields the following

\begin{theorem}\label{Thm5.4}
Let $S$ be a locally noetherian $\BF_q$-scheme and let $b,c:S\to C$ be two $\BF_q$-morphisms. Let $\ul\CE=(\CE,\wt\CE,j,\tau,b,c)$ be a Drinfeld-Anderson shtuka on $C$ of rank $nr'$ and dimension $d'$ over $S$. Then the contravariant functor
\begin{eqnarray*}
\ul T:\es\Sch_{/(C'\times C')\times_{(C\times C)}S} & \longto & \Sets \\[2mm]
\Bigl(S',f:S'\to S \qquad\es\qquad\mbox{ }& \mapsto & \Bigl\{\,\text{Isomorphism classes of pairs $(\ul\CE',\alpha)$ where} \\[-2mm]
(b',c'):S'\to C'\times_{\BF_q} C'\;\Bigr)& & \quad\bullet\es \ul\CE' \text{ is a Drinfeld-Anderson shtuka of rank $r'$, }\\[-2mm]
& & \quad\es\,\text{ dimension $d'$, pole $b'$, and zero $c'$ over $S'$ and}\\[1mm]
& & \quad\bullet\es \alpha:f^\ast\ul\CE\isoto\pi_\ast\ul\CE' \text{ is a fixed isomorphism}\,\Bigr\}
\end{eqnarray*}
is representable by a finite $S$-scheme. \qed
\end{theorem}

The following results are proved analogously to Theorem~\ref{Thm3.6} and Lemma~\ref{Lemma3.7}.

\begin{theorem}\label{Thm5.5}
The restriction of coefficients morphism $\pi_\ast:\DASp^{r',d'}\to\DAS^{nr',d'}$ satisfies the valuative criterion for properness.\qed
\end{theorem}

\begin{lemma}\label{Lemma5.6}
Let $R$ be a valuation ring with fraction field $K$ and let $f:\Spec K\to\Spec R$ be the induced morphism. Let $\ul\CE$ and $\ul\CE'$ be two Drinfeld-Anderson shtuka on $C$ over $\Spec R$ of rank $r$, dimension $d$, pole $b:\Spec R\to C$, and zero $c:\Spec R\to C$. Let $\alpha:f^\ast\ul\CE\to f^\ast\ul\CE'$ be an isomorphism over $\Spec K$. Then $\alpha=f^\ast\beta$ for a unique isomorphism $\beta:\ul\CE\isoto\ul\CE'$ over $\Spec R$. \qed
\end{lemma}

\noindent
{\it Remark.}
Again these results say in the language of stacks that the restriction of coefficients 1-morphism $\pi_\ast:\DrShtp^{r',d'}\to\DrSht^{nr',d'}$ is proper but in general not representable.


\vfill

\parbox[t]{8.2cm}{ 
Urs Hartl  \\ 
Mathematisches Institut  \\ 
Albert-Ludwigs-Universit{\"a}t Freiburg  \\ 
Eckerstr.\ 1  \\ 
D -- 79104 Freiburg\\ 
Germany  \\[0.1cm] 
e-mail: urs.hartl@math.uni-freiburg.de 
} 
\parbox[t]{6.7cm}{
Markus Hendler  \\ 
Mathematisches Institut  \\ 
Albert-Ludwigs-Universit{\"a}t Freiburg  \\ 
Eckerstr.\ 1  \\ 
D -- 79104 Freiburg\\ 
Germany  \\[0.1cm] 

}

\end{document}